\documentclass[12pt,a4paper]{amsart}
\usepackage{latexsym}
\usepackage{amsfonts}
\usepackage[mathscr]{eucal}
\usepackage{epsfig}
\usepackage{a4wide}
\usepackage{ifthen}

\newboolean{draft}
\setboolean{draft}{false}

\def\@copyright{}
\newtheorem{theorem}{Theorem}[section]

\newtheorem{lemma}{Lemma}[section]

\newtheorem{proposition}{Proposition}[section]

\newtheorem{remark}{Remark}[section]

\newcommand{\N}{ \mathbb{N} }

\newcommand{\R}{ \mathbb{R} }

\newcommand{\trunc}[1]{ {\lfloor #1 \rfloor} }

\newcommand{\calE}{\mathcal{E}}

\newcommand{\calK}{\mathcal{K}}
\newcommand{\calL}{\mathcal{L}}

\newcommand{\calM}{\mathcal{M}}

\newcommand{\eins}{{\mathbf 1}}

\newcommand{\vecnull}{{\mathbf 0}}

\newcommand{\Var}{{\mbox{Var\,}}}

\newcommand{\Cor}{{\mbox{Cor\,}}}

\begin{document}

\pagestyle{plain}
\begin{titlepage}
\pagenumbering{arabic}
\title[???]{}
\date{07/2005}
\end{titlepage}

\maketitle

\begin{center}
  \Large
  RANDOM WALKS WITH DRIFT - \\A SEQUENTIAL APPROACH
\end{center}
\vskip 1cm
\begin{center}
  Ansgar Steland\footnote{Address of correspondence: Ansgar Steland, Ruhr-Universit\"at Bochum,
Fakult\"at f\"ur Mathematik, Mathematik 3 NA 3/71, Universit\"atsstr. 150,
D-44780 Bochum, Germany.}\\
  Fakult\"at f\"ur Mathematik\\
  Ruhr-Universit\"at Bochum, Germany\\
  ansgar.steland@rub.de
\end{center}

\begin{abstract}
In this paper sequential monitoring schemes to detect nonparametric drifts are
studied for the random walk case. The procedure is based on a kernel smoother.
As a by-product we obtain the asymptotics of the Nadaraya-Watson estimator and
its associated sequential partial sum process under non-standard sampling.
The asymptotic behavior differs substantially from the stationary situation,
if there is a unit root (random walk component). To obtain meaningful asymptotic
results we consider local nonparametric alternatives for the drift component.
It turns out that the rate of convergence at which the drift vanishes determines
whether the asymptotic properties of the monitoring procedure are determined by
a deterministic or random function. Further, we provide a theoretical result about 
the optimal kernel for a given alternative.
\\[1cm]
\end{abstract}

\newpage
\section*{Introduction}

Many economic time series are non-stationary, and analysts have to take
account of that fact. A time series can be trend-stationary or have 
a random walk component (difference-stationarity). In the first case
shocks are temporary, whereas shocks to a random walk are permanent.
For unit root tests we refer to Dickey and Fuller (1979), Phillip (1979), 
Phillips and Perron (1988), Bierens (1997), and Breitung (2002).
Often, in particular for financial data, the unit root hypothesis 
can not be rejected, and then we are interested to detect as soon as possible
a change-point
where the time series is affected by an additional deterministic drift term. 
The problem discussed in this article should not be mixed up with the
so-called random walk hypothesis which addresses a different issue,
namely whether future values are predictable using past values. 
For that problem we refer to
French and Roll (1986), Fama and French (1988), Lo and MacKinlay (1988), 
Poterba and Summers (1988), and Jegadeesh (1991). 

An important property of a random walk is that there are stochastic
trends which can be mixed up with deterministic trends. 
Nevertheless, a random walk, i.e., a stochastic
trend can be overlayed by a deterministic trend component. 
Hence we study the problem to detect a nonparametric drift
component in a random walk. We assume that the observations
$ Y_{N,1}, Y_{N,2}, \dots, Y_{N,N} $ arrive sequentially and
$$
  Y_{N,n+1} = Y_{N,n} + m_{N,n} + u_n, \qquad n = 1, \dots, N, \quad N \in \N,
$$
where $ u_n $ are i.i.d. innovations with $ E( u_n ) = 0 $ and $ 0 < \Var( u_n ) < \infty $.
For the weak distributional limits presented in this paper the i.i.d. assumption 
can be relaxed by a weak condition discussed in detail in Section~1, which allows, e.g.,
for correlated time series with GARCH effects.
 To study asymptotic properties analytically, we will model the
deterministic drift $ m_{N,n} $ more explicitly. However, the detection procedure
will not depend on a specification of the alternative, but decides after each new
observation whether to continue with observations or whether to stop and reject
the null hypothesis of no drift. In any case we stop no later than after the $N$th
observation, where $N$ is done in advance.

Whereas a posteriori methods aim at estimating consistently 
the time point where the mean changes and therefore employ data before
and after the change point, sequential monitoring methods use only past and 
current data, aiming at the detection of a change as soon as possible.
The a posteriori approach is well studied. For example,
Kim and Hart (1998) propose a nonparametric approach to test for a change
in a mean function when the data are dependent. Predictive tests for
structural change with unknown changepoint have been studied in Ghysels, Guay and
Hall (1997). The analysis of multiple structural changes in linear models has
been discussed, e.g., in Bai and Perron (1998).
Yakir, Krieger and Pollak (1990) use the data after the change for optimization.
Hu\v{s}kov\'{a} and Slab\'{y} (2001) studied nonparametric multiple change point
detection based on kernel-weighted means similar as studied in this article.
Kernel-weighted averages have also been discussed by
Ferger (1994b, 1994c, 1995, 1996) and Brodsky and Darkhovsky (1993, 2002), 
where the latter
examines a posteori and monitoring procedures. 
Sequential monitoring procedures to control for the derivative of a process mean 
have been studied in Schmid and Steland (2000).
Results for jump-preserving smoothers can be found in Chiu {\sl et al.} (1998),
Pawlak and Rafaj\l owicz (2000, 2001), Rue {\sl et al.} (2002),
Steland (2002c, 2004a, 2005a), and Pawlak, Rafaj\l owicz and Steland (2004).
For the application of $U$-statistics we refer to Ferger (1994a, 1997),
Gombay and Horv\v{a}th (1995), and Horv\'{a}th and Hu\v{s}kov\'{a} (2003).

The contribution of the present paper is to study sequential smoothers to monitor
random walks to detect deterministic drifts, and to contrast the results
to former work about stationary processes (Steland 2004b, 2005b).
Whereas in the stationary case the normed delay of the procedure converges
to a deterministic constant, for a random walk the relevant (kernel-weighted)
partial sums have a different convergence rate. Hence, the statistics have to be
scaled appropriately to obtain well-defined limit distributions.
Further, depending on the rate of convergence of the local alternative, we obtain
a deterministic or stochastic limit under the alternative.
As a by-product, we provide the asymptotic law of the Nadaraya-Watson estimator. 
Our approach via kernel-weighted sequential partial sum processes 
yields asymptotic results for both the {\em classic fixed sample design} 
and the {\em sequential sampling design}. Compared to classic nonparametric regression,
the monitoring framework as suggested by Wald (1947), Siegmund (1985), Brodsky and Darkhovsky (1993), and many others, assumes observations at fixed time points.

The paper is organized as follows. 
Section~\ref{Model} discusses the random walk model
with local drift and the proposed monitoring procedure.
Section~\ref{RWAR1} gives a brief discussion of the asymptotics for a stationary AR(1) process.
Section~\ref{Asymptotics} provides the new results about the control statistic under the random walk
model for both the null hypothesis and the alternative. 
We also discuss general time designs in Section~\ref{RWGeneralTimeDesigns}, which
allow to thin a time series with respect to time.
The results are applied in Section~\ref{Threshold}
to derive the related results for the sequential stopping procedures.
Section~\ref{OptKernel} studies the question of optimal kernel choice.
Finally, in Section~\ref{sims} we study the accuracy of the asymptotic distributions by simulations.

\section{Model, method, and assumptions}
\label{Model}

We aim at detecting a nonparametric trend starting at a so-called
change-point (break-point) in the presence of a pure random walk
without drift. In this section we explain in detail the model,
the proposed method and required assumptions.

\subsection{Non-stationary time series model}

Assume the data $ Y_1, \dots, Y_N $, $ N \in \N $, arrive sequentially,
$$
  Y_{N,n+1} = Y_{N,n} + m_{N,n} + u_n, \qquad 1 \le n \le N, \ N \in \N,
$$
where $ \{ u_n \} $ is a sequence of innovation terms with
$ E( u_n ) = 0 $ and common variance $ 0 < \sigma^2 < \infty $.
We assume that the observation $ Y_n $ is taken at time $ t_n $,
where $ \{ t_n \} $ denotes a deterministic and ordered sequence
of time points. For convenience, we assume $ t_n = n \in \N $.
Generalizations to other designs are straightforward and discussed 
in subsection~\ref{RWGeneralTimeDesigns}.

We will study a detection procedure which does not depend on a
specification of the drift $ m_{N,n} $. The null hypothesis (in-control model)
is that $ m_{N,n} $ vanishes for all $ N, n \in \N $, and in this case
$ Y_{N,n} = Y_n $. The alternative says that starting at a change-point $ t_q $
specified below the mean changes.
Our limit theorems work under the following sequence of 
alternative models (out-of-control models) for the drift term. We assume
\begin{equation}
\label{RWModelM}
  m_{N,n} = m_0( [t_n-t_q]/h_N ) h_N^\beta, \qquad n \in \N,\ h > 0,
\end{equation}
where $ h = h_N $ is a sequence of positive constants with
\begin{equation}
\label{NhLim}
  N/h_N \to \zeta \in [1,\infty), \quad \mbox{as $ N \to \infty $}.
\end{equation}
$ m_0 $, called {\em generic alternative}, is a continuous function
such that $ m_0(t) = 0 $ for $ t \le 0 $ and $ m_0(t) \ge 0 $ for $ t > 0 $.
$ m_0 = 0 $ corresponds to the null hypothesis. 
The function $m_0$ is given by nature and unknown to us. However, in many
applications it may be possible to define, e.g., a worst-case scenario in terms of $m_0$.
Then our results can be used to assess the performance of the procedure
under that scenario. $ \beta \in (-1,0] $ is a tuning parameter which controls the rate of 
convergence.
If $ m_0(t) > 0 $, $ t \in (0,t^*) $, for some $ t^* > 0 $, then 
there is a change at time $ t_q $. $t_q$ is called {\em change-point}.
In this paper we address the following two change-point models.

Change-point model CP1: Having in mind applications where it is reasonable to assume
that a change may occur at a fixed given date, e.g., when a firm publishes its
balance sheet, it is assumed that $ t_q = q \in \N $ is a fixed integer. Consequently,
if $ m_0 $ does not vanish, for each fixed $ N $  there is a change, but the
percentage of pre-change observations tends to $0$, as $ N $ tends to $ \infty $.
It will turn out that in this case the asymptotic limit depends on the 
function $m_0$, but not on the change-point.

Change-point model CP2: This approach, which is well established in the literature, assumes
that the change occurs after a fixed fraction of the data, i.e.,
$$
  t_q = t_{Nq} = \trunc{ N \vartheta }, \quad \text{for some $ \vartheta \in (0,1) $}.
$$
Here and in the sequel we
denote by $ \trunc{x} $ the largest integer less or equal to $ x \in \R $.
Under this model the asymptotic limit will depend on the change-point parameter $ \vartheta $,
too.

\begin{remark}
Let us briefly discuss our approach to define local alternatives nonparametrically
 more precisely. We may
write $ m_{N,n} = \widetilde{m}_{N}(t_n;\beta) $, if
$
  \widetilde{m}_{N}(t;\beta) = m_0( [t-t_q]/h ) h_N^\beta.
$
Hence, since $ N/h \to \zeta $, 
for each fixed $ t $ we have $ \lim_{N \to \infty} \tilde{m}_N(t;\beta) = m_0(0) = 0 $. Provided $ m_0(t) $ is twice differentiable at $ t = 0 $ 
with $ m_0''(0) < \infty $, we have
$
  \widetilde{m}_N(t;\beta) = m_0'(0) (t-t_q) h^{\beta-1} + O( h^{\beta-2} ).
$
Thus, the underlying drift tends to zero at rate $ h^{\beta-1} $, point-wise.
\end{remark}

Although in this article we do not discuss the case of
dependent innovations in detail, our results work under the following
general assumption. 

{\bf Assumption (A):} The stationary sequence $ \{ u_n \} $ ensures that the 
partial sum process $ N^{-1/2} \sum_{i=1}^{\trunc{Nr}} u_i $, $ r \in [0,1] $,
converges weakly to scaled Brownian motion, $ \sigma B(r) $, for some
constant $ 0 < \sigma < \infty $ which is determined by 
$ \sigma^2 = \lim_{N \to \infty} E( N^{-1/2} \sum_{i=1}^{N} u_i )^2 $. 

It is worth to discuss assumption (A).
First, note that it covers weakly dependent innovations as arising
in stationary ARMA or GARCH models, 
provided certain additional conditions are fulfilled. In particular, Basrak, Davis
and Mikosch (2003) have shown that GARCH($p,q$) models,
$ Y_n = \sigma_n \epsilon_n $, $ \sigma_n^2 = \alpha_0 + \sum_{i=1}^p \alpha_i
Y_{n-i}^2 + \sum_{j=1}^q \beta_j \sigma_{n-j}^2 $, where $ \{ \epsilon_n \} $ 
are i.i.d. with $ E \epsilon_n = 0 $ and $ E \epsilon_n^2 = 1 $, 
are strictly stationary and strongly mixing
with geometric rate, if $ \alpha_0 > 0 $, $\sum_i \alpha_i + \sum_j \beta_j < 1 $
and $ E \ln^+ | \epsilon_1 | < \infty $, provided the series is started with its
stationary distribution. For a general sufficient condition for (A) in terms of
the $ \alpha $-mixing coefficients $ \{ \alpha(k) : k \in \N \} $ of 
$ \{ u_n \} $ we refer to Herrndorf (1985), which in particular yields (A) provided 
there exists some $ \delta > 0 $ such that $ E |u_1|^{2+\delta} < \infty $ and
$ \sum_{k=1}^\infty \alpha(k)^{\delta/2+\delta} < \infty $.
Finally, note that
this assumption is often considered as a nonparametric definition
of an $ I(0) $ process (e.g. Davidson, 2002). 

\subsection{The monitoring procedure}

We monitor the time series by a sequential kernel smoother
$$
  \widehat{m}_{n} = \sum_{i=1}^n K_h( t_i - t_n) Y_i \, / \, \sum_{i=1}^n K_h(t_i-t_n)
$$
which employs only past and current data. 
The associated kernel-weighted sequential partial sum process is defined as
\begin{eqnarray*}
  \widehat{m}_{N}(s) &=&   \frac{ \sum_{i=1}^{\trunc{Ns}} K_h( t_i - t_{\trunc{Ns}}) Y_i }
                  { \sum_{i=1}^{\trunc{Ns}} K_h( t_i - t_{\trunc{Ns}}) }, \qquad s \in [0,1].
\end{eqnarray*}
$h$ is a bandwidth parameter given in advance 
and $ K_h(z) = h^{-1}K(z/h) $ the rescaled version
of the smoothing kernel $K$. If $K$ vanishes outside the interval $ [-1,1] $,
$ h $ equals the number of past observations used by the procedure.
To obtain meaningful results, namely weak limits, under alternatives, it
turns out that the smoothing parameter $h$ should converge to $ \infty $,
as $ N \to \infty $, i.e.,
$ h = h(N) \uparrow $. It turns out that $ h(N) $ and the sequence $ h_N $ appearing
in the definition of the local alternatives should satisfy 
$ \lim_{N \to \infty} h(N) / h_N = c $ for some constant $ c > 0 $. That
constant can be absorbed in the unknown function $ m_0 $. Thus, for simplicity
we assume $ h(N) = h_N $. The asymptotic framework is parameterized in the
maximum sample size $ N \to \infty $ under the condition (\ref{NhLim}).

Note that the random function $ \widehat{m}_N( \circ ) $ is an element of the Skorokhod space $ D[0,1] $,
consisting of all right-continuous functions with left-hand limits.
We will denote convergence in distribution of random variables
and random vectors by $ \stackrel{d}{\to} $. 
Weak convergence in the space $ D[0,1] $ will be denoted by $ \Rightarrow $. 

The time series is now monitored by the truncated stopping rule
$$
  S_{N} = \inf\{ 1 \le n \le N : T_n(N) > c \}, \qquad T_n(N) = c(h,N) \widehat{m}_{n},
$$
with $ \inf \emptyset = N $. Here $ c(h,N) $ is a scaling function
to be chosen later, and $ T_n(N) $ is the rescaled sequential smoother. 
Note that $ S_{N} $ is the index of the first time point 
where the kernel smoother exceeds the threshold (critical value) $ c $. The monitoring
procedure is truncated, i.e., we stop monitoring at $ N $.
Note that the definition of $ S_N $ does not depend on any model specification of
the alternative.

Concerning the smoothing kernel we make the following assumption.
\begin{itemize}
\item[(K)] $ K $ is assumed to be a Lipschitz continuous probability density with
 mean $0$ and finite variance. Let $L$ be the Lipschitz constant, i.e.,
 $$
  | K(z_1) - K(z_2) | \le L |z_1-z_2|
 $$
 holds true for all $ z_1, z_2 \in \R $.
\end{itemize}

For results under the alternative we need the following conditions.

\begin{itemize}
  \item[(M)] $ m_0 $ is a piecewise continuous funtion.
  \item[(KM)] For the function
  $$
    I(x) = \int_0^x K(s-x) \int_0^s m_0(r) \, dr \, ds,
  $$
  assume $ |I(x)| < \infty $ for all $ x \ge 0 $, $ I \in C(\R_0^+) $,
  $ K(\circ) \cdot \int_0^{\circ} m_0 $ has bounded variation, and there exists some 
  $ x^* > 0 $
  such that $ I( x^* ) > c $.
\end{itemize}

{\bf A nuisance-free procedure.} 
It will turn out that the limiting distribution of $ \widehat{m}_{n} $ 
depends on the nuisance parameter $ \sigma^2 $. A simple candidate is
the naive estimator
\begin{equation}
\label{VarEst}
  \widehat{\sigma}_n^2 = \frac{1}{n-1} \sum_{i=2}^n \Delta Y_i^2,
\end{equation}
where $ \Delta Y_i = Y_i - Y_{i-1} $, $ i = 2, \dots, n $.
Recall that $ \widehat{\sigma}_n^2 $ is consistent for $ \sigma^2 $ under the null hypothesis, if
$ \{ \Delta Y_n \} $ is a linear process, $ \Delta Y_n = \sum_{j=-\infty}^{\infty} \psi_j Z_{n-j} $
where $ \{ Z_j \} $ are i.i.d(0,$\eta^2$) with $ E Z_j^4 < \infty $ and
coefficients satisfying $ \sum_{j=-\infty}^\infty | \psi_j | < \infty $
(Brockwell and Davis, 1991, Prop. 7.3.4). 

A better choice may be Gasser's estimator which is based on a local linear fitting 
procedure (Gasser {\sl et al.}, 1986.) 
Define the pseudo-residuals
$$
  \widetilde{\varepsilon}_n = 0.5 \Delta Y_{n-1} + 0.5 \Delta Y_{n+1} - \Delta Y_n
$$
and note that $ E \widetilde{\varepsilon}_n^2 = ED^2(n,h) + (3/2) \sigma^2 $, where
$
  D(n,h) = (1/2) (m_{n-1} - m_n + m_{n+1} - m_n).
$
By (\ref{RWModelM}) $ D^2(n,h) = O( h^{2 \beta -2} ) $, if $ h \to \infty $, provided $ m_0 $
is twice continuously differentiable. This yields the following proposition.

\begin{proposition} If $ m_0 $ is twice continuously differentiable, the estimator
$$
  \widetilde{\sigma}_{n}^2 = \frac{2}{3(n-2)} \sum_{i=2}^{n-1} \widetilde{\varepsilon}_i^2
$$
is asymptotically unbiased, as $ h \to \infty $, if $ \{ u_n \} $ are i.i.d. with
existing second moment. 
\end{proposition}

Thus, whereas the estimator $ \widehat{\sigma}_n^2 $
tends to overestimate the variance, $ \widetilde{\sigma}_{n}^2 $ may produce more reliable estimates. A related estimator is Rice's estimator given by
$ 1/(2[n-2]) \sum_{i=2}^{n-1} \Delta^2 Y_{i+1} $.

Thus, we may use the asymptotically nuisance-free control statistic
$
  T_{n}^*(N) = s_n^{-1} T_{n}(N).
$
where $ s_n $ is one of estimators discussed above.

\section{Asymptotics for stationary AR processes}
\label{RWAR1}

Before turning our attention to the random walk case, let us briefly discuss the
situation for a stationary $ AR(1) $ process. The asymptotic behavior follows from
general results obtained for stationary $ \alpha $-mixing sequences of innovations,
but the resulting formulas are slightly different and less explicit.

In this section we assume that $ Y_{N,1}, \dots, Y_{N,N} $ are observations
arriving sequentially and
$$
  Y_{N,n+1} = a Y_{N,n} + m_{N,n} + u_n, \qquad n = 1, \dots, N, \quad N \in \N,
$$
where the AR parameter $ a $ satisfies $ |a| < 1 $,
$ \{ u_n \} $ is an i.i.d. sequence of innovations with $ E(u_n) = 0 $ and
$ 0 < \Var(u_n) = \sigma^2 < \infty $. The deterministic drift component is given by
$$
  m_{N,n} = m_0( [t_n-t_q]/h ),
$$
with $ m_0 $ as in the introduction, but at this point we put $ \beta = 0 $. 
Note that $ Y_n = Y_{N,n}, \ n \in \N, $ is stationary under $ H_0 $.

We have
$
  Y_{n+1} = \sum_{i=0}^\infty a^i m_{n-i} + \sum_{i=0}^\infty a^i u_{n-i},
$
where $ \sum_{i=0}^\infty a^i u_{n-i} $ is a stationary process
with autocovariance function
$$
  r_0(k) = \sigma^2 a^k /(1-a^2), \qquad |k| \in \N_0,
$$
thus being $ \alpha $-mixing with geometric rate. 

Under the null hypothesis $ H_0: m_0 = \vecnull $ we may apply Theorem~3.1 of
Steland (2004b) to obtain weak convergence at the usual rate $ N^{1/2} $, i.e.,
\begin{equation}
\label{ARH0}
  \frac{h}{N^{1/2}} \widehat{m}_N(s) 
  \Rightarrow \mathbb{M}(s), \qquad \mbox{in $D[0,1]$},
\end{equation}
as $ N \to \infty $, where $ \mathbb{M}_\zeta(s) $ is a centered Gaussian process with 
correlation kernel given by
$$
  \Cor( \mathbb{M}_\zeta(s), \mathbb{M}_\zeta(t) ) = C_\zeta(s,t)\ /\ 
  \left( \zeta^2 \int_0^{\zeta s} K(z-\zeta s) \, dz \int_0^{\zeta t} K(z - \zeta t) \, dz \right),
$$ 
for $ 0 \le s \le t \le 1 $, with
$$
  C_\zeta(s,t) = \lim_{N \to \infty} \sum_{i=1}^{\trunc{Ns}} \sum_{j=1}^{\trunc{Nt}} K_h(t_i-t_{\trunc{Ns}}) K_h(t_j - t_{\trunc{Ns}}) 
  \frac{ \sigma^2 a^{|i-j|} }{1-a^2}.
$$
Due to the Lipschitz continuity of $ K $, the sample paths of $ \mathbb{M}_\zeta $ 
are continuous w.p. $1$. Note that
$$
  I(m_0) = \lim_{N \to \infty} \sum_{j=0}^N a^j m_0( [n-j]/h ) < \infty,
$$
if $ \int m_0^2(s) \, ds < \infty $.
Now a similar argument as in Theorem~3.3 of Steland (2004b) shows that
under the alternative the process in (\ref{ARH0}) diverges at the rate $ N^{1/2} $,
since
\begin{eqnarray*}
  h N^{-1} \widehat{m}_N(s) 
  & = & \frac{h}{N} \sum_{i=1}^{\trunc{Ns}} K_h(t_i-t_{\trunc{Ns}}) \sum_{j=0}^i a^j m_0( [n-j]/h ) + o_P(1) \\
  & = & O \left( I(m_0) \int_0^{\zeta s} K(z - \zeta s) \, dz \right).
\end{eqnarray*}
These results have also immediate implications for the sequential stopping rules. 
If 
$$
  \mu_\zeta(s) = P-\lim_{N \to \infty} h N^{-1} \widehat{m}_N(s),
$$
it can be shown that for any fixed $ 0 < \kappa < 1 $
$$
  N^{-1} \inf\{ \trunc{\kappa N} \le n \le N : h N^{-1/2} \widehat{m}_n > c \}
  \stackrel{P}{\to} \inf \{ \kappa \le s \le 1 : \mu_{\zeta}(s) > c \},
$$
as $ N \to \infty $, i.e., the normed delay converges to a deterministic quantity.

\section{Asymptotics for random walks}
\label{Asymptotics}

Now we study the asymptotic behavior of the Nadaraya-Watson estimator
$ \widehat{m}_{n} $ under the random walk model as introduced in
Section~\ref{Model}.
Note that our asymptotic framework differs
from the usual framework in nonparametric regression. We do not
assume that the time points $ \{ t_i \} $ get dense in any finite time interval or
are distributed according to a density, which ensures
that we may let the bandwidth $h$ tend to $ 0 $ at a certain rate.
Instead we assume a fixed time design taking account of the fact
that time series are commonly observed at a fixed time scale.
Thus, as a by-product we provide the asymptotic laws of the
Nadaraya-Watson type smoothing under the sampling design of the present paper.
We formulate the results for equidistant observations, i.e., $ t_n = n $, and
discuss more general time designs in Section~\ref{RWGeneralTimeDesigns}.

The results of this section about the Nadaraya-Watson process
$ \widehat{m}_N(s) $, $ s \in [0,1]$, are preparations for the analysis
of the stopping time $ S_N $, but since they are interesting in their own right
we discuss them in detail here. In particular, the interesting relationship between the
(qualitative) asymptotic behavior and the convergence rate of the local
alternative are properties of that underlying process.

\subsection{Limit theory under the null hypothesis}

We first study the asymptotic distributions under the null hypothesis that
we deal with a random walk without drift. The limit distributions are 
centered Gaussian processes and centered normal distributions, respectively.

\begin{theorem} 
\label{Th1} Assume (A) and (K).
Under the null hypothesis $ H_0: m_0 = \vecnull $ we have
$$
  h N^{-3/2} \widehat{m}_{N} \stackrel{d}{\to} 
    \frac{\sigma \int_0^1 K( \zeta (r-1) ) B(r) \, dr}
         { \zeta \int_0^1 K( \zeta(r-1) ) \, dr },
$$
as $ N \to \infty $. The associated partial sum process converges weakly
$$
  h N^{-3/2} \widehat{m}_{N}(s) 
  \Rightarrow \calM_\zeta(s) = \frac{ \sigma \int_0^s K( \zeta(r-s) ) B(r) \, dr } 
                   { \zeta \int_0^s K( \zeta(r-s) ) \, dr }, \quad \text{in $D[0,1]$},
$$
as $ N \to \infty $. The limit process is continuous w.p. $ 1 $.
\end{theorem}

Observe that for $ \sigma = 1 $ the limit process $ \calM_\zeta(s) $ is distributed according
to a $ N(0,\sigma_K^2) $ distribution with variance given by
$$
  \sigma_K^2(s) = \frac{
    \int_0^1 K(\zeta(s-1)) \biggl[ \int_0^s t K( \zeta (t-1) ) \, dt +
                                s \int_s^1 K( \zeta (t-1) ) \, dt \biggr] \, ds 
    }
    {
      \left( \zeta \int_0^s K( \zeta(r-s) ) \, dr \right)^2
    }
$$
which can be calculated explicitly for any given kernel (Shorack and Wellner (1986), p. 42). 
The following table provides some
values of $ \sigma_M^2 = \sigma_K^2(1) $ for the Gaussian kernel, the Epanechnikov kernel given
by $ K_{Epan}(z) = (3/4)(1-z^2), $ for $ z \in [-1,1] $, and the (standardized) Laplace kernel,
which is defined by $ K_{Lap}(z) = (1/\sqrt{2})e^{-\sqrt{2}|z|}$, $ \ z \in \R $.

\begin{table}[h]
\begin{center}
\begin{tabular}{lccccccc}
  Kernel           &   \multicolumn{7}{c}{$\zeta$} \\ 
                   &  $10$   &   $5$     &  $ 4 $   &  $2$  &   $1.5$  &   $1.2$  &   $ 1 $ \\ \hline
  Gaussian         &  $0.0089$&  $0.0310$&  $0.0449$& $0.1242$& $0.1913$ & $0.2754$&  $0.3775$ \\
  Laplace          &  $0.0089$&  $0.0316$&  $0.0463$& $0.1443$& $0.2310$ & $0.3353$&  $0.4578$ \\
  Epanechnikov     &  $0.0095$&  $0.0359$&  $0.0545$& $0.1857$& $0.2921$ & $0.3968$&  $0.4857$ \\ \hline \\
\end{tabular}
\caption{Asymptotic variances for several choices of the kernel and $ \zeta = \lim N/h $.}
\label{RWTab1}
\end{center}
\end{table}

Theorem~\ref{Th1} suggests the following confidence interval
\begin{equation}
\label{RWConfInt}
  \widehat{m}_{N} \pm z_{1-\alpha/2} \sigma_K h^{-1} N^{3/2}
\end{equation}
which has asymptotic coverage $ 1-\alpha $ under $ H_0 $.
It can be used to perform a preliminary level $ \alpha $ test given data $ Y_1, \dots, Y_N $
before establishing a monitoring procedure. The accuracy of that procedure is studied to some
extent in Section~6.
However, comparing $ \widehat{m}_{Nh} $ with the confidence limits
$ z_{1-\alpha/2} \sigma_K h^{-1} N^{3/2} $ does not ensure well-defined
statistical properties of the associated stopping rule in terms of
the average run length or the normed delay.

\begin{remark} Note that the event $ T_n(N) = h N^{3/2} \widehat{m}_n > c $ stands for
a false alarm at the $n$th time point, if $ m_0 = \vecnull $. It is straightforward to show
$$
  P( h N^{3/2} \widehat{m}_n > c ) = O( h^{-2} N^{3/2} ) = O(N^{-1/2}),
$$
i.e., in our framework the point-wise false-alarm rate tends to $0$, as $ N \to \infty $.
\end{remark}

\subsection{Limit theory under local drifts}

We will now investigate the asymptotic behavior under the (local) alternative model
and both model specifications for the change-point.
It turns out that the result depends qualitatively
on the rate parameter $ \beta $ of the alternative. If $ \beta = -1/2 $, i.e., the
alternative converges at the rate $ h^{-3/2} $ to the null model, we obtain a non-degenerate
Gaussian limit with drift for the process $ h N^{-3/2} \widehat{m}_N(s) $ studied in 
Theorem~\ref{Th1} under the null hypothesis. That process has a proper asymptotic null
distribution.
For a slowly converging alternative ($ \beta = 0 $) corresponding to the rate $ h^{-1} $, 
we have to change the scaling function to obtain a limit. In this case
we obtain stochastic convergence to a non-stochastic function. 
That function determines the asymptotic detection
properties of the proposed procedure. We formulate the results for the partial
sum processes $ \widehat{m}_{N}(s) $, putting $ s = 1 $ yields the
asymptotic laws of the Nadaraya-Watson estimator.

\begin{theorem}
\label{Th3} Assume (A), (K), (M), and (KM). Fix $ 0 < a < 1 $.
Under the alternative $ H_1: m_0 \ge^* 0 $ the following assertions hold true.
\begin{itemize}
\item[(i)] If $ \beta = -1/2 $, we have weakly in $ D[a,1] $,
  $$
  \frac{h}{N^{3/2}} \widehat{m}_{N}(s) \Rightarrow
  \sigma \frac{ \int_0^s K(\zeta(r-s)) B(r) \, dr }{ \zeta \int_0^s K( \zeta(r-1) ) \, dr }
  +
  \frac{  \int_0^s K(\zeta(r-s)) \int_0^{ \zeta r} m_0( t - \zeta \vartheta 1_{CP2}) \, dt \, dr }
       { \zeta^{3/2} \int_0^s K( \zeta(r-s) ) \, dr },
$$
as $ N \to \infty $. Here, $ 1_{CP2} = 0 $ if change-point model CP1 holds, and 
$ 1_{CP2} = 1 $ under model CP2.
\item[(ii)] If $ \beta = 0 $, then 
  $$
    \frac{h^{1/2}}{N^{3/2}} \widehat{m}_{N}(s)
    \stackrel{P}{\to} \frac{ \int_0^s K( \zeta(r-s) ) \int_0^{\zeta r} m_0( t - \zeta \vartheta 1_{CP2} ) \, dt \, dr }
    { \zeta^{3/2} \int_0^s K( \zeta(r-s) ) \, dr },
  $$
as $ N \to \infty $. Again, $ 1_{CP2} = 0 $ if change-point model CP1 holds, and 
$ 1_{CP2} = 1 $ under model CP2.
\end{itemize}
\end{theorem}

\begin{remark} Note that the asymptotic limit depends on the change-point parameter
$ \vartheta $ if model CP2 holds. Under model CP1 the limit is free of $ t_q $, which
is a consequence of $ t_{q}/h = o(1) $ and continuity of $ m_0 $.
\end{remark}

\begin{remark} It is worth noting that procedures based on the partial sum process
  $ \widehat{m}_{N}(s) $ are able to detect a drift if the function
  $$
    \mu_\zeta(s) =  \int_0^s K( \zeta(r-s) ) \int_0^{\zeta r} m_0( t - \zeta \vartheta 1_{CP2} ) \, dt \, dr
  $$
  is positive for some interval of $s$-values.
\end{remark}

\begin{remark} Note that (ii) implies that the statistic $ h N^{-3/2} \widehat{m}_N $ diverges
  under local alternatives corresponding to $ \beta = 0 $ at the rate $ h^{1/2} $.
\end{remark}

\section{General time designs}
\label{RWGeneralTimeDesigns}

Let us briefly discuss more general time designs for the choice of the time 
points $ t_n $.
In some applications the following monitoring approach may be possible and reasonable. 
We monitor the process at equidistant time points $1, 2, \dots $ 
until either the procedure provides a signal, or we have reached the time horizon
(maximum sample size) $ N $. Here we assume that the time unit is chosen appropriately,
e.g., one day or one week. Intuitively, to detect a change as soon as possible
it should be better to use more recent observations $ Y_i $, i.e. with $ t_n - t_i $ small,
than past observations where $ t_n-t_i $ is large. 
To some extent, this is achieved by the smoothing kernel, which downweights past
data, but a real thinning of the data can only be achieved by an appropriate
selection resp. design of the time points. This means, at the current time $ t_{n,n} = n $
one chooses past time points $ 0 < t_{n,1}, \dots, t_{n,n-1} < t_{n,n} $ where observations are taken. This allows to start with monthly observations and use daily observations at the end of the (current) sample. 
We consider two different approaches corresponding to the two change-point models CP1 and
CP2.

\subsection{Generalized time designs for the CP1 model} 

Assume that
\begin{equation}
\label{RWGeneralTimeDesignRule1}
  t_{n,i} = n F_T^{-1}( i/n ), \qquad i = 1, \dots, n, \ n \in \N,
\end{equation}
where $ F_T $ is a continuously differentiable d.f. with support $[0,1] $. 
Clearly, if $ F_T $ is the d.f. of the uniform distribution on $ [0,1] $, we obtain $ t_{n,i} = i $.
Nonlinear choices of $ F_T $ allow to ensure that
past or more recent observations dominate the sample. Note that
$ F_T^{-1} $ defines a sampling scheme which is rolled over the time axis:
At each time $n$ the time points $ t_{n,1}, \dots, t_{n,n-1} $ are chosen
according to the scheme (\ref{RWGeneralTimeDesignRule1}).

Under model CP1, a Taylor expansion yields $ t_{nq} = (F_T^{-1})'(0) q + o(1) $
provided $ F_T^{-1} $ is continuously differentiable.
Thus, if $ (F_T^{-1})'(0) > 0 $, the underlying (asymptotic) change-point equals
$ (F_T^{-1})'(0) t_q $, whereas for $ (F_T^{-1})'(0) = 0 $ the sequence of
change-points vanishes asymptotically, i.e., the detection problem is made easier
as $n$ increases.

The associated Nadaraya-Watson process is given by
$$
  \widehat{m}_N(s) = \frac{ \sum_{i=1}^{ \trunc{Ns} } K_h( t_{ \trunc{Ns}, i} - \trunc{Ns} ) Y_i }
                          { \sum_{i=1}^{ \trunc{Ns} } K_h( t_{ \trunc{Ns}, i} - \trunc{Ns} ) },
  \qquad s \in [0,1],
$$
where again $ \trunc{Ns} $ plays the role of the current time point.
It is straightforward to check that the proofs of Theorem~\ref{Th1} and Theorem~\ref{Th3} still work. Now the limit process under the null hypothesis is given by
$$
   \calM_{\zeta,F_T}(s)  = \frac{ \sigma \int_0^s K( \zeta s[ F_T^{-1}(r/s) - 1] ) B(r) \, dr }
        { \int_0^s K( \zeta s[ F_T^{-1}(r/s) - 1] ) }, \qquad s \in [0,1].
$$
The drift term appearing in Theorem~\ref{Th3} changes to
$$
  \mu_{\zeta,F_T}(s) = \frac{ \int_0^s K( \zeta s[ F_T^{-1}(r/s) - 1 ] ) \int_0^r m_0(t) \, dt \, dr }
       { \zeta^{3/2} \int_0^s K( \zeta s [ F_T^{-1}(r/s) - 1 ] ) \, dr }, \qquad s \in [0,1].
$$

\begin{remark} \label{Rem}
In practice, it may be necessary to
use the time point $ t_{n,j}^* \in \{ t_{n,1}^*, \dots, t_{n,m}^* \} $ nearest to $ t_{n,i} $,
where $ \{ t_{n,j}^* \} $ denotes the finest discrete time scale available. Then, 
(\ref{RWGeneralTimeDesignRule1}) defines a selection rule for the time points $ \{ t_{n,j}^* \} $.
\end{remark}

\subsection{Generalized time designs for the CP2 model}

It is easy to see that the generalized time
design above makes not much sense under model CP2. One may consider the
following modification, which is easier to apply, but lacks the authentic idea
to allow for schemes which use more observations near each current time $ n $. Assume
\begin{equation}
\label{RWGeneralTimeDesignRule2}
  t_{N,i} = N F_T^{-1}( i/N ), \qquad i = 1, \dots, N, 
\end{equation}
where $ F_T $ is a continuously differentiable d.f. with support $ [0,1] $. 
Here, given the maximum sample size $N$,
the time design scheme is set up only once, i.e, the selected time points
do not change with the current time $n$. Since under model CP2 the
change-point is given by $ t_q = t_{Nq} = \trunc{N \vartheta} $, we obtain
$$ 
  t_{Nq} = N F_T^{-1}( \trunc{N \vartheta} / N ) 
$$ 
yielding
$ t_{Nq} / N \to F_T^{-1}( \vartheta ) $. This means, the (asymptotic) change-point 
parameter is transformed by $ F_T^{-1} $, and it appears in the asymptotic limit.
The associated Nadaraya-Watson process is now defined by
$$
  \widehat{m}_N(s) 
    = \frac{ \sum_{i=1}^{ \trunc{Ns} } 
        K_h( N/h [F_T^{-1}( i/N ) - F_T^{-1}(\trunc{Ns}/N)] ) Y_i 
      }
      { \sum_{i=1}^{ \trunc{Ns} } K_h( N/h [F_T^{-1}( i/N ) - F_T^{-1}(\trunc{Ns}/N)] ) }
  \qquad s \in [0,1],
$$

A straightforward calculation shows that the drift term now changes to
$$
  \mu_{\zeta,F_T}(s) 
   = \frac{ \int_0^s K( \zeta s[ F_T^{-1}(r/\zeta) - F_T^{-1}(s) ] ) 
     \int_0^{\zeta r} m_0( \zeta[ F_T^{-1}( t/\zeta ) - F_T^{-1}( \vartheta/\zeta) ] ) \, dt \, dr }
       { \zeta^{3/2} \int_0^s K( \zeta s [ F_T^{-1}(r/\zeta) - F_T^{-1}(s) ] ) \, dr }, \qquad s \in [0,1].
$$
Note that Remark~\ref{Rem} also applies to the time design scheme 
(\ref{RWGeneralTimeDesignRule2}).

\section{Sequential detection rules}
\label{Threshold}

Let us now discuss the implications of the results of Section~\ref{Asymptotics} for the
stopping rule
$
  S_N = \inf \{ 0 \le n \le N : T_{n} > c \}.
$
Note that $ S_N $ can be written in terms of the sequential partial sum processes. Indeed,
$
  S_N = N \cdot \inf \{ 0 \le s \le 1 : c(h,N) \widehat{m}_N(s) > c \}.
$
For asymptotic results under local alternatives we also consider the
stopping rule
$$
  S_N^{(a)} = \inf \{ \trunc{Na} \le n \le N : c(h,N) \widehat{m}_N( n/N ) > c \}
$$
where $ a \in (0,1) $ is a fixed constant.
Again notice that $ S_N^{(a)} $ can be written as
$ N \cdot \inf \{ a \le s \le 1 : c(h,N) \widehat{m}_N(s) > c \}. $

\subsection{Limit theory under the null hypothesis}

The following theorem provides the null distribution of the stopping rules.

\begin{theorem} 
\label{Th3a}
Assume (A), (K), and $ H_0: m_0 = \vecnull $ (random walk without drift).
\begin{itemize}
\item[(i)]
  If $ T_{n}(s) = c(h,N) \widehat{m}_{n}(s) $  with scaling factor 
  $ 
    c(h,N) = h N^{-3/2}, 
  $
  the normed stopping time $ S_{N} / N $ 
  converges in distribution to the random variable
  $$
    S_\zeta = \inf \biggl\{ s \in [0,1] : 
      \frac{ \sigma \int_0^s K(\zeta(r-s)) B(r) \, dr }
           { \zeta \int_0^s K( \zeta(r-s)) \, dr } 
      > c \biggr\},
  $$
  as $ N \to \infty $.
\item[(ii)] The limiting laws of the nuisance-free versions correspond to the special case $ \sigma = 1 $.
\end{itemize}
\end{theorem}

These results can be used to choose the threshold (critical value) $ c $ from the asymptotic
distribution. For example, we may simulate trajectories from the limiting processes 
and determine for each trajectory
the smallest $ s $ such that the threshold $ c $ is exceeded. This gives an approximation
of the distribution of $ S_{N} $ which can be used to choose $ c $ to ensure that,
e.g., the average run length equals a prespecified value.

\subsection{Limit theory under local drifts}

The following results summarize our findings under local alternatives and
give interesting insights into the asymptotic properties of the procedure.
In particular, we see how the smoothing kernel and the generic alternative
$ m_0 $ jointly affect the performance of the procedures.

\ifthenelse{\boolean{draft}}{
We need the following uniform convergence result for the expectation of the numerator
of $ h N^{-3/2} \widehat{m}_N(s) $, namely,
$$
  \mu_N(s) = E\biggl( h N^{-3/2} \sum_{i=1}^{\trunc{Ns}} K_h(t_i-t_{\trunc{Ns}}) Y_i \biggr).
$$

\begin{lemma}
  We have
  $$
    \sup_{a \le s \le 1}
    \left| h^{-1/2} \mu_N(s) - \int_0^s K( \zeta(r-s) ) \int_0^r m_0( t - \zeta \vartheta 1_{CP2}) \, dt \, dr
    \right| \to 0,
  $$
  as $ N \to \infty $.
\end{lemma}
}{}

\begin{theorem} 
\label{Th4}
Assume (A), (K), (M), and (KM) (random walk with local drift).
Fix $ a \in (0,1) $.
\begin{itemize}
\item[(i)] 
  Suppose $ \beta = -1/2 $. If $ T_{n}(s) = c(h,N) \widehat{m}_{N}(s) $ with scaling factor
  $
    c(h,N) = h N^{-3/2},
  $
  the normed stopping time $ S_{N}^{(a)} / N $
  converges weakly to the random variable
  $$
    S(\zeta) = \inf \{ s \in [a,1] : W_\zeta(s) > c \},
  $$
  where the stochastic process $ W_\zeta(s) $ is given by
  \begin{eqnarray*}
    W_\zeta(s) &=& \frac{ \sigma \int_0^s K( \zeta(r-s) ) B(r) \, dr }{ \zeta \int_0^s K( \zeta(r-s) ) \, dr } 
                  + 
       \frac{ \int_0^s K(\zeta(r-s)) \int_0^r m_0(t - \zeta \vartheta 1_{CP2}) \, dt \, dr }
            { \zeta^{3/2} \int_0^s K( \zeta(r-s) ) \, dr },
  \end{eqnarray*}
  as $ N \to \infty $.
  \item[(ii)]
  Suppose $ \beta = 0 $. If $ T_{n}(s) = c(h,N) \widehat{m}_{N}(s) $
  with scaling factor 
  $
    c(h,N) = h^{1/2} N^{-3/2},
  $ 
  the normed stopping time $ S_{N}^{(a)}/N $ 
  converges in probability to the non-stochastic asymptotic normed delay
  \begin{eqnarray*}
    S^*(\zeta;K;m_0) &=&
    \inf \left\{ s \in [a,1] : 
        \frac{ \int_0^s K(\zeta(r-s)) \int_0^r m_0(t - \zeta \vartheta 1_{CP2}) \, dt \, dr }
             { \zeta^{3/2} \int_0^s K( \zeta(r-s) ) \, dr }
          > c \right 
    \},
  \end{eqnarray*}
  as $ N  \to \infty $.
\end{itemize}
\end{theorem}

This theorem says that the stopping rule relying on the control statistic $ h N^{-3/2} \widehat{m}_N $, 
which has a proper limit under $ H_0 $, has a nondegenerate limit distribution under local alternatives
converging to $0$ at the rate $ h^{-3/2} $. If, however, we consider alternatives with rate
$ h^{-1} $, which is the appropriate rate in the stationary case (see Steland, 2004b), and change the scaling function, we obtain a deterministic limit
$ S^*(\zeta;K;m_0) $, the asymptotic normed delay, as in the case of a stationary process.

\section{Optimal kernel choice}
\label{OptKernel}

Suppose the critical value $ c $ is a fixed constant chosen by the data analyst. For example,
when analyzing a time series representing financial risk measured in terms of a
currency unit, $ c $ may be a psychological price. 
Then $ S_N $ stands for the time point where that price is reached for the first time.

Assuming the change point model CP1, Theorem~\ref{Th4}~(ii) motivates to 
examine whether optimal kernels exist which
minimize the asymptotic normed delay $ S^*(\zeta;K;m_0) $ for a given 
alternative $ m_0 $ representing a worst case scenario.
Recall that this deterministic quantity appears as the limit if the alternative model
converges to the null model at the rate $ h^{-1} $, whereas for the faster rate
$ h^{-3/2} $ we obtained a stochastic limit. From a practical viewpoint considering
the conditions for a slower convergence to $0$ may provide a better approximation
to reality.

First note that for a finite set of candidate kernels, $ \{ K_1, \dots, K_M \} $, 
we can simply plot the $ M $ corresponding curves
$$
  y_l(s) = \frac{ \int_0^s K_l(\zeta(r-s)) \int_0^r m_0(t) \, dt \, dr }
                { \zeta \int_0^s K( \zeta (s-r) ) \, dr },
  \quad l = 1, \dots, M,
$$
and use the kernel which provides the smallest $ s $ where the critical value $ c $ is exceeded.
For the case of detecting a drift in a stationary process Steland (2002a) provides a real data analysis of
credit risk data, where this simple procedure yields a detection rule which signals 
the change one time point earlier. For a Bayesian view on the problem of kernel optimization see Steland (2002b).

Although we can provide a solution to the problem of optimal kernel choice, the results seem to be 
of limited practical use, since we can identify the optimal kernel only for a finite 
interval around $ 0 $. 
Nevertheless, from a theoretical point of view it is interesting to know that
both the asymptotic normed delay and the optimal kernel can be calculated explicitly for any given generic alternative $ m_0 $.

Let $ \calK $ denote a class of probability densities with expectation $ 0 $,
which is uniformly Lipschitz continuous, i.e.,
$$
  \sup_{K \in \calK} | K(z_1) - K(z_2) | \le L |z_1-z_2|, \qquad \forall z_1, z_2 \in \R,
$$
holds for some constant $ L > 0 $. The problem is to find a kernel $ K^* \in \calK $
such that the corresponding asymptotic normed delay, $ S^*(\zeta;K^*;m_0) $, satisfies
$$
  S^*( \zeta; K^*; m_0 ) = \inf \{ S^*(\zeta; K; m_0) : K \in \calK \}.
$$
Such a pair $ (K^*, S^*(\zeta;K^*;m_0)) $ is called {\em optimal}.
Using optimization techniques presented in detail in Steland (2004b), one can
establish  the following theorem which provides a way to calculate the optimal asymptotic normed delay and provides the optimal kernel $ K^* $.

\begin{theorem} 
\label{OptimalKernel}
Suppose that for all $ s \in [0,1] $
$$
  0 < \int_0^s \left( \int_0^r m_0(t) \, dt \right)^2 \, dr < \infty.
$$
\begin{itemize}
\item[(i)] The optimal asymptotic normed delay is given by
  $$
    S^*(\zeta,K^*;m_0) = \inf \left\{ s \in [0,1] :
    \frac{ \int_0^{s} \left( \int_0^{r} m_0(t) \, dt \right)^2 \, dr }
         { \int_0^s \int_0^r m_0(t) \, dt \, dr } 
    > c \right\}.
  $$
\item[(ii)]
The optimal kernel $ K^* $ satisfies
$$
  K^*(z) = \frac{ \int_0^{z/\zeta + S^*(\zeta;K^*;m_0)} m_0(t) \, dt }
                { 2 \int_0^{\infty} \int_0^r m_0(t) \, dt \, dr }
$$
for arguments $ z \in [- \zeta S^*(\zeta;K^*;m_0), \zeta S^*(\zeta;K^*;m_0)] $.
\end{itemize}
\end{theorem}

\section{Simulations}
\label{sims}

To study the accuracy of the asymptotic distributions 
of the detection procedures, we simulated random walks, $ \{ Y_n \} $, where
$ Y_0 = 0 $, and $ Y_{n+1} = Y_n + u_n $ with $ \{ u_n \} $ i.i.d. $ N(0,\sigma^2) $,
$ \sigma = 1 $.
To estimate the nuisance parameter $ \sigma^2 $ 
we assumed that an additional prerun random walk of length $ h = 10 $ was given.

Figure~1 shows $20$ realizations of the kernel-weighted sequential partial sum process,
$ \widehat{m}_{N}(s) $, $ s \in [0,1] $, for $ N = 100 $ and $ h = 50 $ and its asymptotic approximation via the kernel-weighted
integral over Brownian motion using $ \zeta = 2 $. 
The sequential detection procedure $ S_N $ can be visualized by drawing a horizontal line (control limit) 
at $ c $.
The first intersection of the process and the control limit is the run length.

To study the accuracy of the asymptotic null distribution we performed simulations to assess the coverage
of the confidence interval based on $ \widehat{m}_N $ and average run lengths (ARL) of 
the stopping rule $ S_N $. We focus on the ARL, since it may the most
common criterion to design monitoring procedures for practical applications. Note, however, that
our results also allow to design procedures which control the type I error rate.

Table~\ref{RWSim1} reports the simulated coverage probabilities of the confidence interval 
defined in (\ref{RWConfInt}) for a Gaussian kernel and a nominal 
coverage of $ 0.95 $ under the null hypothesis. The results for the Epanechnikov and Laplace kernel, respectively, 
were in close agreement and are not reported here. Each value is estimated by 10.000 repetitions.
The asymptotic 
variance is estimated using the estimator (\ref{VarEst}) and $ \sigma_K^2 $ as given in 
Table~\ref{RWTab1}. It can be seen that even for $ h << N $ and small $N$ coverage is good.

\begin{table}
\begin{tabular}{cccccc} \hline
$\zeta$   &     10 &     50 &   100  &  250   &  500   \\ \hline
10 & 0.9502 & 0.9496 & 0.9523 & 0.9502 & 0.9471  \\
5  & 0.9481 & 0.9514 & 0.9478 & 0.9489 & 0.9534  \\
4  & 0.9475 & 0.9525 & 0.9474 & 0.9473 & 0.9515  \\
2  & 0.9408 & 0.9468 & 0.9480 & 0.9458 & 0.9518  \\
1.5& 0.9350 & 0.9431 & 0.9516 & 0.9512 & 0.9485  \\
1.2& 0.9320 & 0.9453 & 0.9523 & 0.9477 & 0.9518  \\
1  & 0.9301 & 0.9526 & 0.9470 & 0.9494 & 0.9504 \\ \hline \\
\end{tabular}
\caption{Coverage probabilities of a $0.95$-confidence interval for random walks.}
\label{RWSim1}
\end{table}

\begin{figure}[t]
\includegraphics[width = 12cm, height = 10cm]{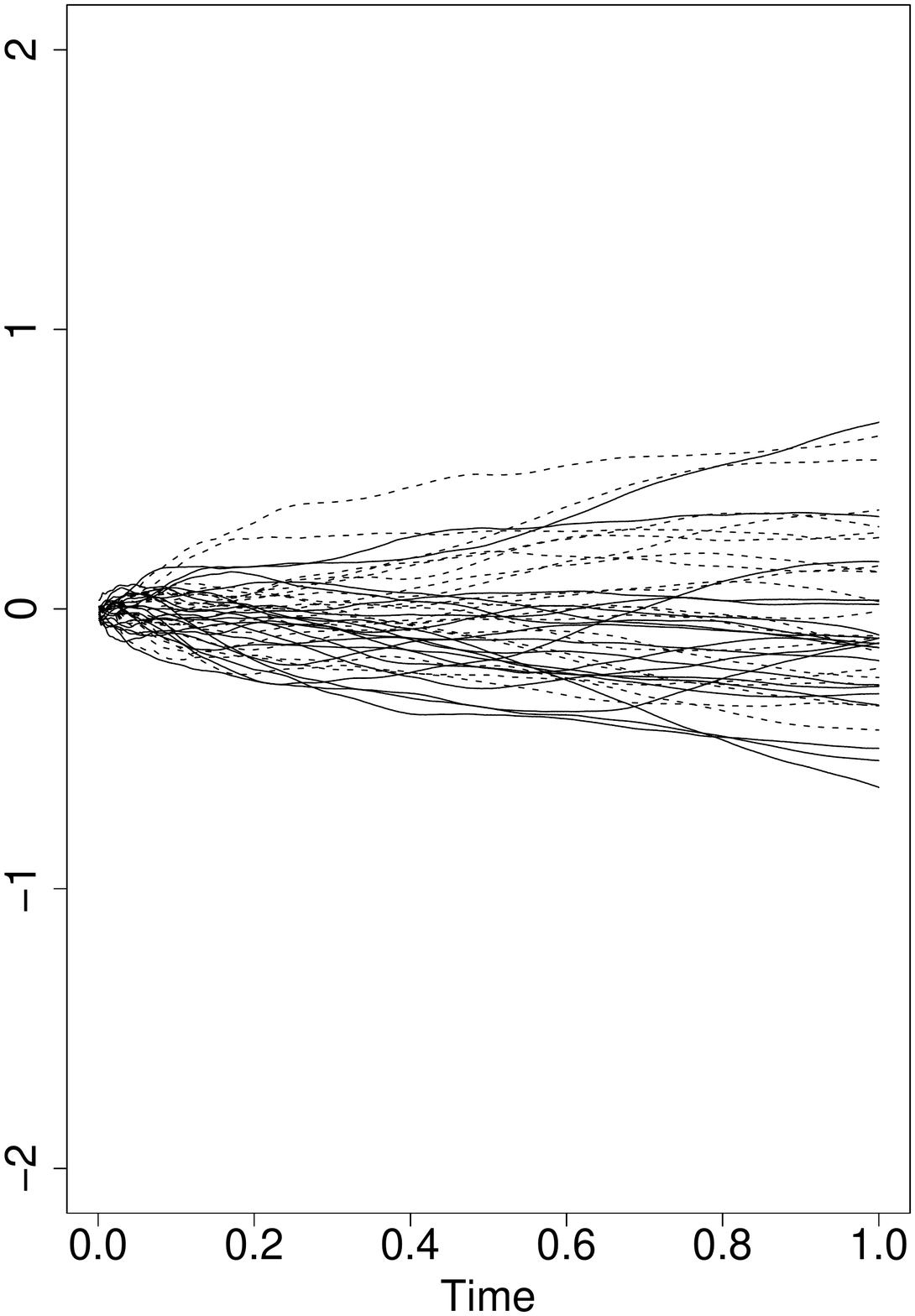}
\caption{$20$ realizations of the kernel-weighted sequential process $ \widehat{m}_{nh}(s),\ s \in [0,1] $,
  (bold line) and its asymptotic limit (dashed line).}
\label{Fig1}
\end{figure}

\begin{figure}[t]
\includegraphics[width = 12cm, height = 10cm]{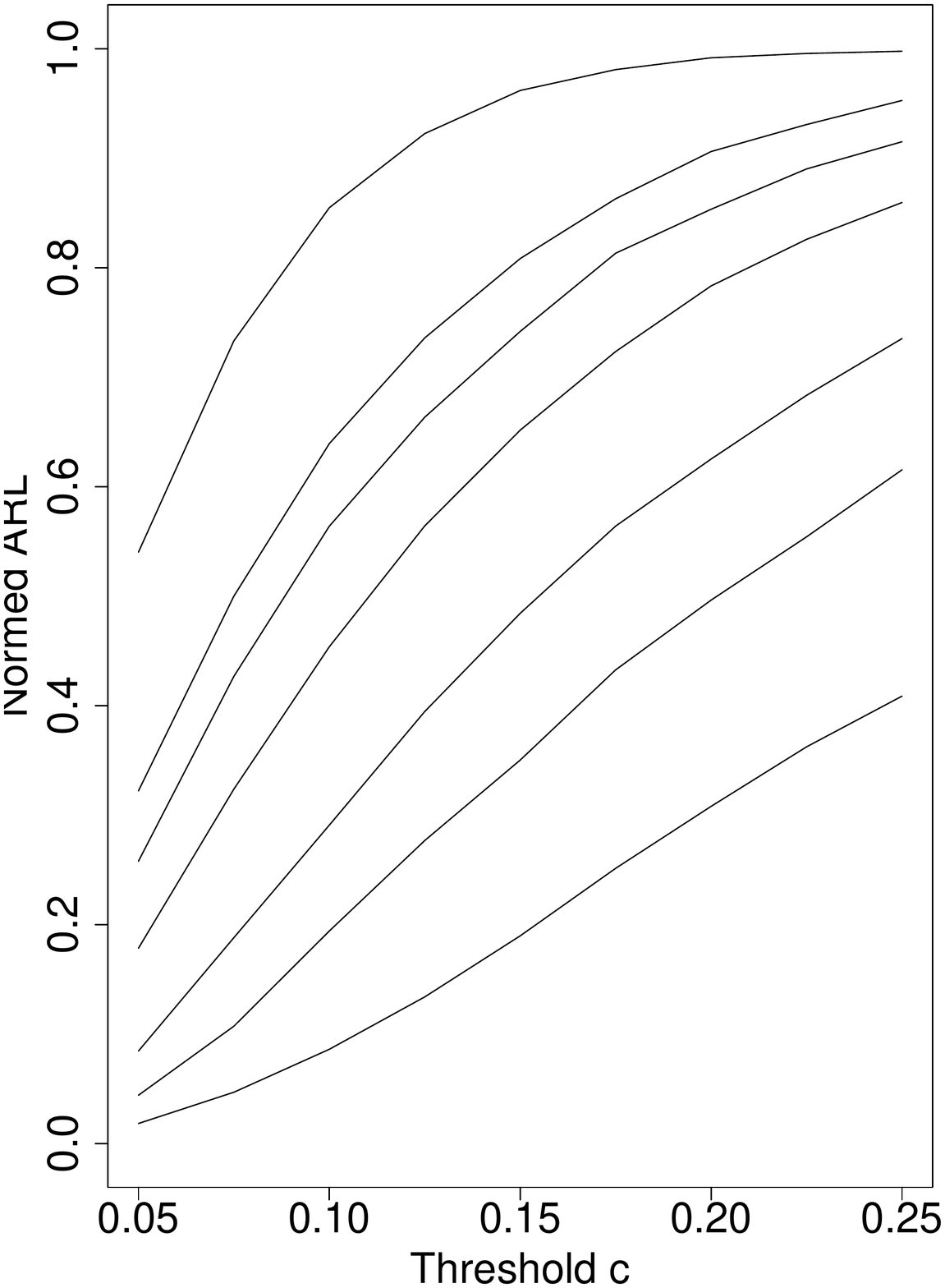}
\caption{Normed ARL curves for $ h N^{-3/2} \widehat{m}_{nh} $ using the Gaussian kernel. 
  $ \zeta $ attains the values $ 1 $ (bottom curve), 1.5, 2, 3, 4, 5, and 10 (top curve).}
\label{Fig3}
\end{figure}

In order to simplify the application of the proposed sequential monitoring procedure we provide
curves to obtain approximate critical values to achieve a prespecified ARL,
$ E_0( S_N ) $, under the null hypothesis $ H_0: m_0 = \vecnull $.
Figure~\ref{Fig3} provides curves of the normed ARL $ a_0 = E_0(S_N) / N $ as a function of $ c $, 
i.e., $ a_0 = a_0(c) $.
For given $ (N,h) $ use the curve for $ \zeta \approx N/h $ and
determine $ c $ graphically with $ N c \approx a_0(c) $.

How accurate is that approximation? To gain some insight we compared
the asymptotic distribution of the stopping time 
$$
  S_\zeta = \inf \{ 0 \le s \le 1 : \int_0^s K( \zeta(r-s) ) B(r) \, dr  / \int_0^s K( \zeta(r-s) ) \, dr > c \}
$$
with the true distribution of the normed stopping time
$$
  S_{Nh} / N = \inf \{ 1 \le n \le N : \widehat{\sigma}_n^{-1} h N^{-3/2} \widehat{m}_{nh}(s) > c \} / N
$$
in terms of the ARL. Each ARL was approximated using $ 10,000 $ trajectories.

\begin{figure}[t]
\includegraphics[width = 6cm, height = 7cm]{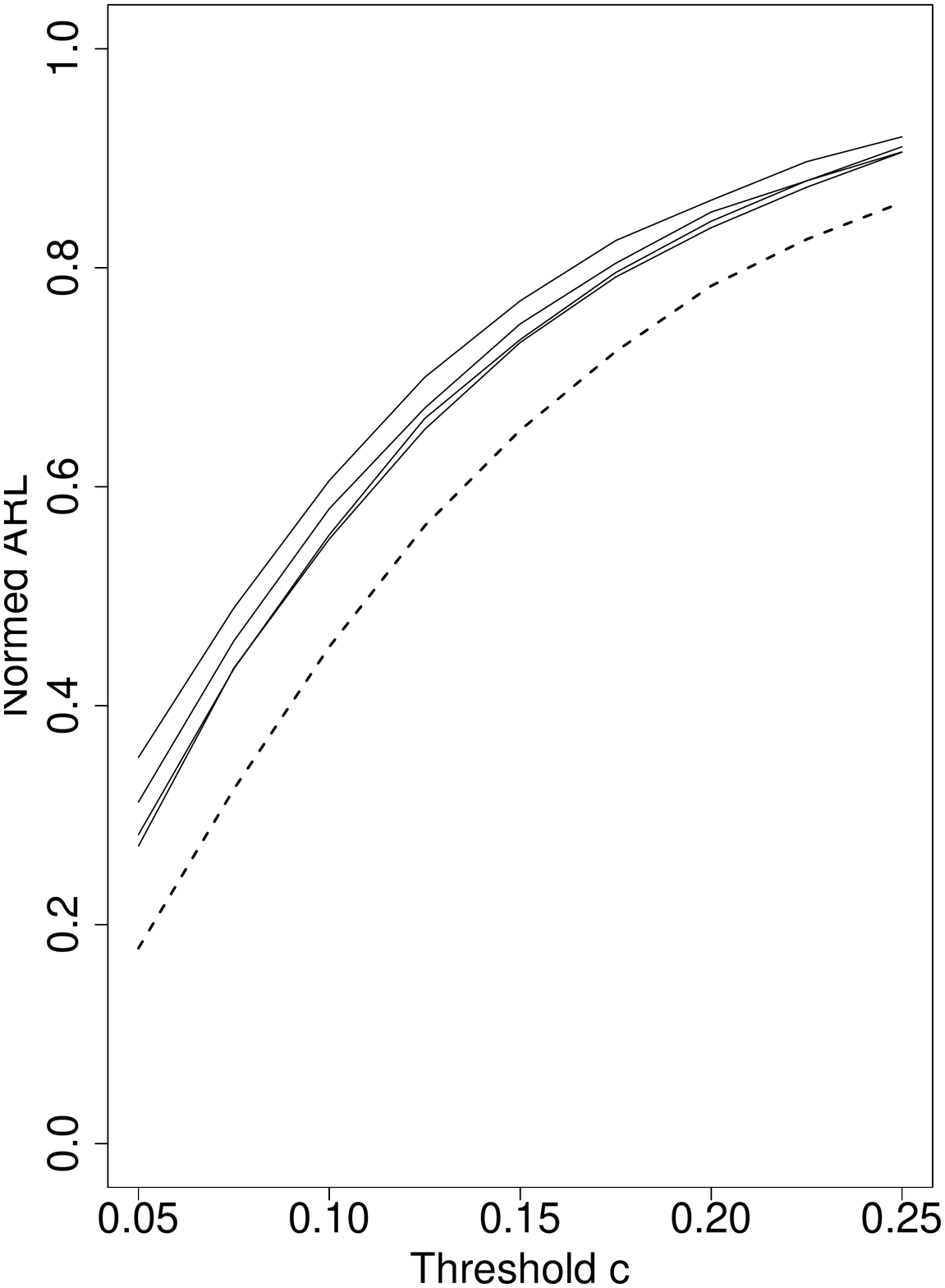}
\includegraphics[width = 6cm, height = 7cm]{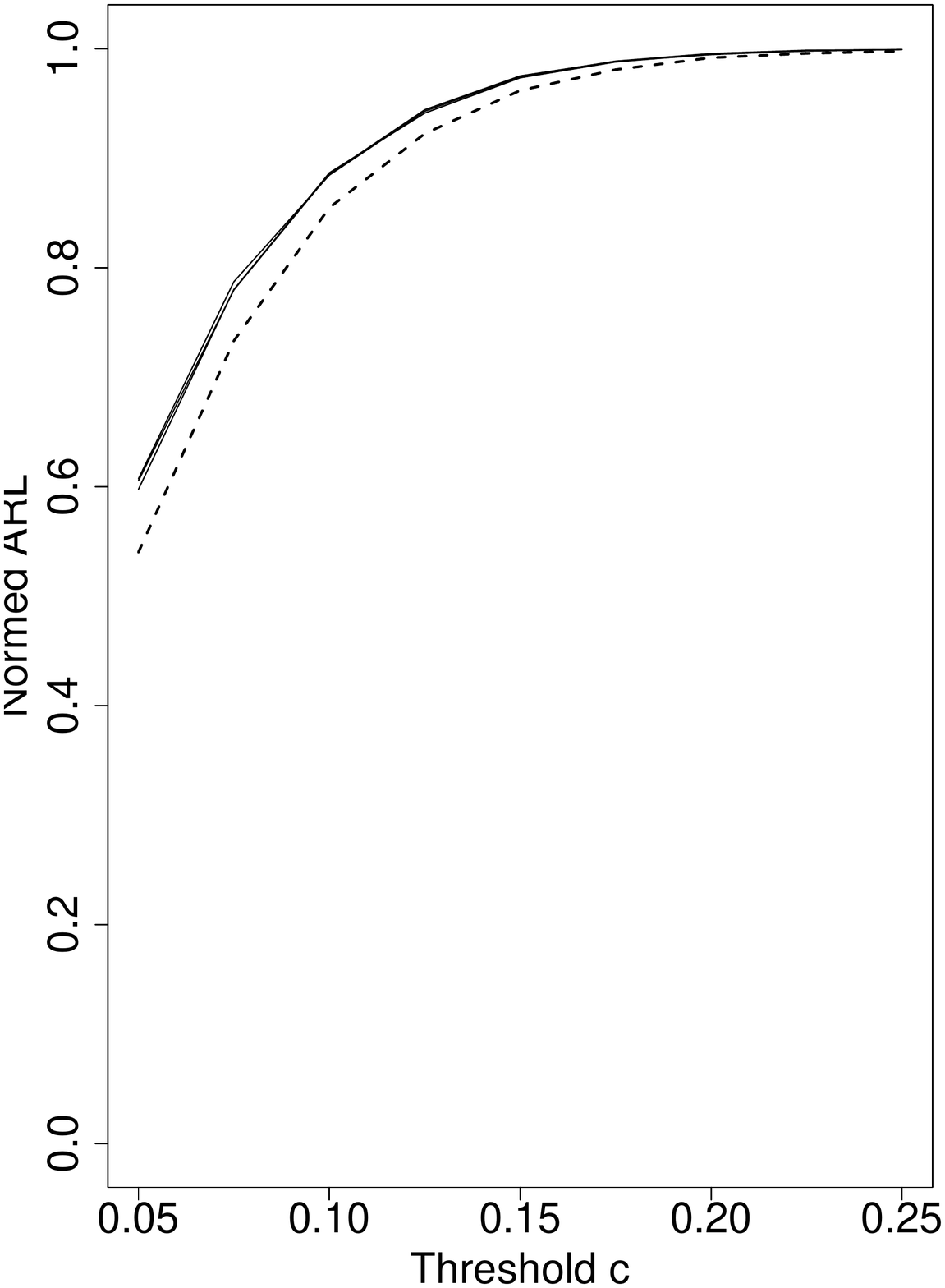}
\caption{Normed ARL curves for the nuisance-free control statistic 
  $ T_{nh}^* = \widehat{\sigma}_n^{-1} hN^{-3/2} \widehat{m}_{nh} $ using the Gaussian kernel. 
  $ h $ takes on the values $ 10, 20, 50, 100 $, $ N = \zeta h $. 
  {\em Left panel:} $ \zeta = 3 $. {\em Right panel:} $ \zeta = 10 $.
  The dashed curves represents normed ARLs of the asymptotic distribution.}
 \label{Fig4}
\end{figure}

Figure~\ref{Fig4} provides the results. For $ h \in \{ 10, 20, 50, 100 \} $, $ N = \zeta h $, 
and $ \zeta = 3 $ (left panel) and $ \zeta = 10 $ (right panel) the corresponding normed-ARL curves
are shown. It can be seen that the curve representing the asymptotic critical values are below the
simulated true curves. This means, the asymptotic critical values yield conservative procedures.
The accuracy seems to be better for large values of $ \zeta $, i.e., if $ h $ is small compared to $ N $.

\section*{Acknowledgements}

The author is grateful to an anonymous referee and an associate editor for valuable suggestions which improved the paper. The support of the
Deutsche Forschungsgemeinschaft (DFG) (SFB 475, {\em Reduction of Complexity in Multivariate Data Structures}) is acknowledged.


\section*{APPENDIX: PROOFS}

\ifthenelse{\boolean{draft}}{
Let us briefly discuss the notion of weak convergence which is used in this article
to formulate our result on the asymptotic distributional limits of the relevant
random functions.
Weak convergence of a sequence $ \{ X_n(s) : s \in [0,1] \} $ of stochastic processes, i.e.,
elements of the space $ D[0,1] $, to
a process $ X(s) $, also an element of $ D[0,1] $ is 
denoted by $ X_n(s) \Rightarrow X(s) $ and understood in the sense that
for all continuous and bounded test functions $ h: D[0,1] \to \R $, 
$Eh(X_n) \to Eh(X) $, as $ n \to \infty $. This definition is equivalent to
convergence in distribution of all finite-dimensional distributions, i.e., for
all distinct time-points $ t_1, \dots, t_k \in [0,1] $ the sequence of random vectors
$ ( X_n(t_1), \dots, X_n(t_k) ) $ converges in distribution to $ (X(t_1),\dots, X(t_k)) $,
as $ n \to \infty $, and for all $ \varepsilon > 0 $ there exists a compact set 
$ B \subset D[0,1] $ such that all $ X_n $ are contained in $ B $ with probability 
$ 1 - \varepsilon $. For more details we refer to Billingsley (1968), Pollard (1985) and
Hall and Heyde (1980).
}
{
In this paper we work with weak convergence (denoted by $ \Rightarrow $) of elements of the space $ (D[0,1],d) $ where $d$ 
is the Skorokhod metric. For treatments of the general theory we refer to Billingsley (1968), Pollard (1985),
and Vaart and Wellner (1996).
}

\begin{myproof}{of Theorem~\ref{Th1}}
Put $ Y_0 = 0 $ and define
$$
  X_{N}(r;s) = N^{-1} Y_{\trunc{Nr}} K_h( t_{\trunc{Nr}} - t_{\trunc{Ns}} ), \qquad r,s \in [0,1].
$$
Note that $ X_{N}(r;s) $ is a constant on the intervals $ [ \frac{i}{N}, \frac{i+1}{N} ) $ with
value $ N^{-1} Y_i  K_h( t_i - t_{\trunc{Ns}} ) $, $ i = 1, \dots, N $. Therefore, the area under the curve
$ X_{N}(r;s), \ r \in [0,s] $, is given by
$$
  \int_0^s X_{N}(r;s) \, dr = \frac{1}{N^2} \sum_{i=1}^{\trunc{Ns}} K_h(t_i - t_{\trunc{Ns}}) Y_i.
$$
Using $ Y_{\trunc{Nr}} = \sum_{i=1}^{\trunc{Nr}} u_i $, we have
$  h N^{1/2} X_{N}(r;s) = \frac{1}{\sqrt{N}} \sum_{i=1}^{\trunc{Nr}} u_i \cdot 
    K \left( \frac{ t_{\trunc{Nr}} - t_{\trunc{Ns}} }{ h } \right)$.
Since by assumption (A)
the partial sum process $ N^{-1/2} \sum_{i=1}^{\trunc{Nr}} u_i $ converges weakly
to scaled Brownian motion $ \sigma B(r) $, we may apply the a.s. representation theorem of Skorokhod and
Dudley (Pollard (1984), p. 71) which ensures that there exist versions of the random elements
which converge a.s. in the supnorm. This implies
$$
  \left\| \frac{1}{\sqrt{N}} \sum_{i=1}^{\trunc{N \circ_1}} u_i
    K \left[ \frac{ \trunc{N\circ_1} - \trunc{N \circ_2} }{ h } \right]
  - \sigma B( \circ_1 ) K( \zeta[ \circ_1 - \circ_2 ] )
  \right\|_{D([0,1]\times[0,1])} \to 0,
$$
which proves weak convergence 
$ h N^{1/2} X_{N}(r;s) \Rightarrow \sigma K(\zeta(r-s))B(r) $ in 
$ D([0,1]\times[0,1]) $. By continuity of $ K $, the process
$ \sigma K( \zeta(r-s) ) B(r) $, $ (s,r) \in [0,1] \times [0,1] $,
has continuous and bounded sample paths w.p. $ 1$. Consider
the integral operator $ I $ which maps
an element $ f \in D([0,1]\times[0,1]) $ to the element $ I(f) \in D[0,1] $ given by
$
  I(f)(s) = \int_0^s f(r,s) \, dr,
$
$ s \in [0,1] $.
If $ (f_n) \subset D([0,1] \times [0,1]) $ is a convergent sequence with limit
$ f \in C( [0,1] \times [0,1] ) $, i.e., 
$ d(f_n,f) \to 0 $, as $ n \to \infty $, then we also have $ \| f_n - f \|_\infty \to 0 $,
$ n \to \infty $, yielding $ \| I(f_n) - I(f) \|_\infty \to 0 $, as $n \to \infty $,
i.e., continuity of $ I $. Hence, the continuous mapping theorem yields
\begin{eqnarray*}
  \frac{h}{N^{3/2}} \sum_{i=1}^{\trunc{N \circ_1}} K_h(t_i - t_{\trunc{N \circ}}) Y_i 
    & = & I( h N^{1/2} X_{N}(\circ_2;\circ_1) )( \circ_1 ) \\
    & \Rightarrow &
      \sigma \int_0^{\circ_1} K( \zeta (r-\circ_1) ) B(r) \, dr,
\end{eqnarray*}
weakly in $ D[0,1] $, as $ N \to \infty $. Since additionally,
\begin{equation}
\label{RWSumWeights}
  \sum_{i=1}^{\trunc{Ns}} K_h(t_i - t_{\trunc{Ns}}) 
  \to 
   \int_0^{\zeta s} K( r - \zeta s) \, dr
  = \zeta \int_0^s K( \zeta(r-s) ) \, dr,
\end{equation}
as $ N \to \infty $, the assertions follow.
\end{myproof}


\begin{myproof}{of Theorem~\ref{Th3}}
  A random walk with non-vanishing drift, $ Y_{n+1} = Y_n + m_{nh} + u_n $, can be 
  decomposed as 
  $
    Y_n = \widetilde{Y}_n + \sum_{s=1}^{n-1} m_{sh}, \ n \in \N,
  $
  where $ \widetilde{Y}_n = \sum_{s=1}^{n-1} u_s $ 
  is a random walk based on the innovations $ u_n $ without drift. Hence,
  $$
    h N^{-3/2} \sum_{i=1}^{\trunc{Ns}} K_h( t_i - t_{\trunc{Ns}} ) Y_i
  $$
  can be decomposed as
  $$
    h N^{-3/2} \sum_{i=1}^{\trunc{Ns}} K_h( t_i - t_{\trunc{Ns}} ) \widetilde{Y}_i
      + h N^{-3/2} \sum_{i=1}^{\trunc{Ns}} K_h( t_i - t_{\trunc{Ns}} ) \sum_{j=1}^{i-1} m_{jh}
  $$
  For the first term one may argue as in the proof of Theorem~\ref{Th1} to verify that
  \begin{equation}
  \label{RWH0Conv1}
    h N^{-3/2} \sum_{i=1}^{\trunc{Ns}} K_h(t_i-t_{\trunc{Ns}}) \widetilde{Y}_i
    \Rightarrow
    \sigma \int_0^s K(\zeta(r-s)) B(r) \, dr,
  \end{equation}
  as $ N \to \infty $.
  \ifthenelse{\boolean{draft}}{
    Let us now consider the second term. 
    Observe that $ \{ t_i/h : i = 1, \dots, \trunc{Ns} \} $ forms an
    equidistant partition of $ [1/h, [Ns]/h] $ with size $ 1/h $, 
    where $ [1/h, \trunc{Ns}/h] \to (0,\zeta s) $, as $ N/h \to \zeta $.
    Further, the drift is modelled by
    $$
      m_{n} = m_0( [t_i-t_q]/h ) h^\beta.
    $$
    Therefore, if $ \beta = -1/2 $, we have
  }
  {
    Further, since $ m_n = m_0([t_i-t_q]/h) h^\beta $, $ \beta = -1/2 $ implies 
  }
  \begin{eqnarray*}
    \mu_{N}(s) & = & \frac{h}{N^{3/2}} \sum_{i=1}^{\trunc{Ns}} K_h(t_i - t_{\trunc{Ns}}) \sum_{j=1}^{i-1} m_0( [t_j - t_q]/h ) h^\beta \\
\ifthenelse{\boolean{draft}}{
    & = &
      \frac{ h^{\beta+2} }{N^{3/2}} 
      \frac{1}{h^2} 
      \sum_{i=1}^{\trunc{Ns}} \sum_{j=1}^{\trunc{Ns}} 
        K\biggl(\frac{t_i}{h}-\frac{\trunc{Ns}}{h}\biggr) 
        \eins\biggl( \frac{t_j}{h} \le \frac{t_{i-1}}{h} \biggr) 
        m_0\biggl( \frac{t_j}{h} - \frac{t_q}{h} \biggr) \\
    & \to &
     \frac{1}{\zeta^{3/2}} 
     \int_0^{\zeta s} K(r_1 - \zeta s) \int_0^{r_1} m_0( r_2 ) \, dr_2 \, dr_1 \\
    & = &
      \frac{1}{\zeta^{1/2}} \int_0^s K( \zeta(r-s) ) \int_0^r m_0(t) \, dt \, dr,
}{
    & \to &
      \frac{1}{\zeta^{1/2}} \int_0^s K( \zeta(r-s) ) 
        \int_0^{\zeta r} m_0(t - \zeta \vartheta 1_{CP2}) \, dt \, dr,
}
  \end{eqnarray*}
  as $ N \to \infty $, by (K) and (KM) uniformly in $ s \in [a,1] $ (cf. Steland 2004b, Th. 3.3 (ii)).
  Combining this fact with (\ref{RWH0Conv1}) and (\ref{RWSumWeights}) yields
  $$
    \frac{h}{N^{3/2}} \widehat{m}_{nh}(s) 
    \Rightarrow
      \sigma \frac{ \int_0^s K( \zeta(r-s) ) B(r) \, dr }{ \zeta \int_0^s K( \zeta(r-s) ) \, dr }
      + \frac{ \zeta^{-1/2} \int_0^s K( \zeta(r-s) ) \int_0^{\zeta r} m_0(t- \zeta \vartheta 1_{CP2}) \, dt \, dt }
             { \zeta \int_0^s K( \zeta(r-s) ) \, dr },
  $$
  in $ D[a,1] $, as $ N \to \infty $.
  In contrast, if $ \beta = 0 $ we obtain convergence to a
  deterministic quantity, if we change the scaling factor from $ h N^{-3/2} $ to 
  $ h^{1/2} N^{-3/2} $. Indeed, in this case we have
  $$ 
    \frac{h^{1/2}}{N^{3/2}} \sum_{i=1}^{\trunc{Ns}} K_h(t_i-t_{\trunc{Ns}}) \widetilde{Y}_i = o_P(1),
  $$
  uniformly in $ s \in [a,1] $, and for the centering term
  \begin{eqnarray*}
    h^{-1/2} \mu_{N}(s) 
    & = &
    \frac{h^{3/2}}{N^{3/2}} h^{-2}\sum_{i=1}^{\trunc{Ns}} K([t_i-t_{\trunc{Ns}}]/h) \sum_{j=1}^{i-1} m_0([t_j-t_q]/h) \\
    & \to & \zeta^{-1/2} 
    \int_0^s K( \zeta(r-s) ) \int_0^{\zeta r} m_0(t- \zeta \vartheta 1_{CP2}) \, dt \, dr,
  \end{eqnarray*}
  yielding
  $$
    \frac{h^{1/2}}{N^{3/2}} \widehat{m}_{N}(s)
    \stackrel{P}{\to} 
    \frac{ \zeta^{-1/2} \int_0^s K( \zeta(r-s) ) \int_0^{\zeta r} m_0(t- \zeta \vartheta 1_{CP2}) \, dt \, dr }{ \zeta \int_0^s K( \zeta(r-s) ) \, dr },
  $$
  uniformly in $ s \in [a,1] $, as $ N \to \infty $.
\end{myproof}

\begin{myproof}{of Theorem~\ref{Th3a} and \ref{Th4}} 
We verify Theorem~\ref{Th3a} (i), i.e., assuming $ \beta = -1/2 $
and $ c(h,N) = h N^{-3/2} $.
The other assertions are shown along these lines.
Fix $ 0 < a < 1 $. By Theorem~\ref{Th3}~(i) the process 
 $
   c(h,N) \widehat{m}_{N}(s)
 $
 converges weakly in $ D[a,1] $ to the non-stationary
 and a.s. continuous process
 $$
   W_\zeta(s) = \frac{ \sigma \int_0^s K(\zeta(r-s)) B(r) \, dr }{ \zeta \int_0^s K( \zeta(r-s) ) \, dr } 
   + 
   \frac{ \zeta^{-1/2} \int_0^{\zeta s} K(r-\zeta s) \int_0^r m_0(t) \, dt \, dr }
        { \zeta \int_0^s K( \zeta(r-s) ) \, dr },
 $$
 as $ N \to \infty $. Define the functional
 $ \varphi_a : D[0,1] \to D[0,1] $,
 $$
   \varphi_a(f) = \inf \{ a \le s \le 1 : f(s) > c \}, \qquad f \in D[a,1].
 $$
 Clearly, $ \varphi_a | \calE_c $ is continuous w.r.t. $ \| \circ \|_\infty $
 and $ d $, where
 $
   \calE_c = \{ f \in C[0,1] : f(x^*) > c \ \mbox{for some $x^*$} \}.
 $
 By (K) and (M) we have $ W_\zeta \in C[a,1] $ w.p. $ 1 $. Thus, since
 $ S_N^{(a)} / N = \varphi_a( c(h,N) \widehat{m}_N( \circ ) ) $,
 the continuous mapping theorem yields
 $$
   S_N^{(a)}/N \Rightarrow \varphi_a( W_\zeta( \circ ) ) 
    = \inf \{ a \le s \le 1 : W_\zeta(s) > c \},
 $$
 as $ N \to \infty $. Notice that
 $$
   \inf \{ a \le s \le 1 : W_\zeta(s) > c \} > x
   \Leftrightarrow
   \sup_{0 \le s \le x} W_\zeta(s)  \le c.
 $$
 By a.s. continuity of $ W_\zeta $, Theorem 2 of Lifshits (1982) ensures
 that
 $
   \nu_x = \calL( \sup_{0 \le s \le x} W_\zeta(s) )
 $
 can have an atom only at the point
 $$
   \gamma_x = \sup_{0 \le t \le x: Var( W_\zeta(t) ) = 0}
   E W_\zeta(t),
 $$
 vanishes on $ (-\infty, \gamma_x) $, and is absolutely continuous
 on $ (\gamma_x, \infty) $. Since $ \Var( W_\zeta(s) ) > 0 $ if
 $ s > 0 $, $ \nu_x $ is absolutely continuous. Therefore,
 we obtain convergence in distribution, i.e.,
 $$
   P( \inf \{ a \le s \le 1 : c(h,N) \widehat{m}_N(s) > c \} \le x)
   \to P( \inf \{ a \le s \le 1: W_\zeta(s) > c \},
 $$
 as $ N \to \infty $, for all $ x $.
\end{myproof}

\begin{myproof}{of Theorem~\ref{OptimalKernel}}
  Using standard arguments of functional optimization theory,
  we see that $ S^*(\zeta;K;m_0) $ is minimized w.r.t. $ K $ if 
  \begin{equation}
  \label{ToOptim} \tau(K) = 
    \int_0^{s^*} K( \zeta(r-s^*) ) \int_0^r m_0(t) \, dt \, dr \ / \ 
    \int_0^{s^*} K( \zeta(r-s^*) ) \, dr
  \end{equation}
  is maximized w.r.t. $ K \in \calK $, where $ s^* = S^*(\zeta, K^*; m_0) $ 
  denotes the optimal asymptotic normed delay
  (c.f. Steland (2004b)).  Clearly, $ \tau(K) \ge 0  $ is less than or equal to
  $$
  \sqrt{ \int_0^{s^*} K(\zeta(r-s^*))^2 \, dr }
         \sqrt{ \int_0^{s^*} \left( \int_0^r m_0(t) \, dt \right)^2 \, dr } 
    \ / \ 
    \int_0^{s^*} K( \zeta(r-s^*) ) \, dr
  $$
  with equality if and only if
  $$
    \frac{ K( \zeta(r-s^*) ) }
         { \int_0^{s^*} K( \zeta(r-s^*) ) \, dr } 
    = \lambda \int_0^r m_0(t) \, dt, \qquad r \in [0, s^*],
  $$
  for some $ \lambda $. Using $ \int_0^\infty K( \zeta(r-s^*) ) \, dr = 1/2 $ gives
  $$
    \lambda^{-1} = 2 \int_0^{\infty} \int_0^r m_0(t) \, dt \, dr  \int_0^{s^*} K( \zeta(r-s^*) ) \, dr,
  $$
  i.e., the optimal (symmetric) kernel $ K^* $ satisfies
  \begin{equation}
  \label{RWReprOptK}
    K^*( \zeta(r-s^*) ) = \frac{ \int_0^r m_0(t) \, dt  }
                               { 2 \int_0^{\infty} \int_0^r m_0(t) \, dt }, \quad r \in [-s^*,s^*].
  \end{equation}
  Consequently, using $ K(\zeta(r-s^*)) = K(\zeta(s^*-r)) $ and 
  substituting $ z = \zeta s^* - \zeta r $ gives the representation in the theorem
  for $ z \in [-\zeta s^*, \zeta s^*] $. Plugging in $ K^* $ as given in (\ref{RWReprOptK}) in 
  (\ref{ToOptim})
  \ifthenelse{\boolean{draft}}{and noting that
  $$
    \int_0^{s^*} K^*(\zeta(r-s^*))\, dr = \frac{ \int_0^{s^*} \int_0^r m_0(t) \, dr }
                                               { 2 \int_0^\infty \int_0^r m_0(t) \, dr },
  $$
  }{}
  yields immediately
  $$
    \tau(K^*) = \frac{ \int_0^{s^*} K^*(\zeta(s^*-r)) \int_0^r m_0(t) \, dt \, dr }{ \int_0^{s^*} K^*(\zeta(s^*-r)) \, dr} \\
       = \frac{ \int_0^{s^*} \bigl( \int_0^{r} m_0(t) \, dt \bigr)^2 \, dr }
              { \int_0^{s^*} \int_0^r m_0(t) \, dt \, dr }.
  $$
  Therefore, the assertion for the optimal asymptotic normed delay follows.
\end{myproof}

\end{document}